\documentclass[12pt,twoside,reqno]{amsart}
\usepackage{amsmath, amsthm, amscd, amsfonts, amssymb, graphicx}
\let
\oldsection
\section
\renewcommand{\section}{\vspace{8pt plus 4pt}\oldsection}
\usepackage{lipsum}
\usepackage[super]{nth}
\usepackage{xcolor}
\usepackage{tikz}
\usepackage{pgfplots}
\usepackage{mathrsfs}
\usepackage{graphics}
\usepackage{graphicx} 
\usepackage{epsfig}
\definecolor{olive}{rgb}{0.3, 0.4, .1}
\definecolor{fore}{RGB}{249,242,215}
\definecolor{back}{RGB}{51,51,51}
\definecolor{title}{RGB}{255,0,90}
\definecolor{dgreen}{rgb}{0.,0.6,0.}
\definecolor{gold}{rgb}{1.,0.84,0.}
\definecolor{JungleGreen}{cmyk}{0.99,0,0.52,0}
\definecolor{BlueGreen}{cmyk}{0.85,0,0.33,0}
\definecolor{RawSienna}{cmyk}{0,0.72,1,0.45}
\definecolor{Magenta}{cmyk}{0,1,0,0}

\theoremstyle{definition}

\theoremstyle{remark}
\newtheorem{rem}{Remark}[section]
\newtheorem{example}{Example}[section]

\numberwithin{equation}{section}
\begin{document}
\begin{center}\large{{\bf{Some basic properties of G-Calculus and its applications in  numerical analysis}}} 
\vspace{0.5cm}

Khirod Boruah and Bipan Hazarika$^{\ast}$ 

\vspace{0.5cm}
Department of Mathematics, Rajiv Gandhi University, Rono Hills, Doimukh-791112, Arunachal Pradesh, India\\

Email: khirodb10@gmail.com; bh\_rgu@yahoo.co.in
\thanks{$^{\ast}$ The corresponding author.}
\end{center}
\title{}
\author{}
\thanks{{\today}}
\begin{abstract} Objective of this paper is to introduce a new type of calculus which will be called G-Calculus based on non-Newtonian calculus introduced by Grossman and Katz \cite{GrossmanKatz}. The basic difference between geometric calculus defined by Grossman and Katz and the present G-calculus is that Grossman took the values of the argument as $x, x+ h, x+2h,...$ but here in G-calculus we take the values as $x, x\oplus h, x\oplus e^2\odot h, x\oplus e^3\odot h....$ This calculus will have great deal with numerical analysis which are discussed in the last section of this paper.

\parindent=5mm
\noindent{\footnotesize {\bf{Keywords and phrases:}}} Geometric calculus; geometric integers; geometric real numbers.\\
{\footnotesize {\bf{AMS subject classification \textrm{(2000)}:}}} 26A06, 11U10, 08A05, 46A45.

\end{abstract}
\maketitle

\maketitle

\pagestyle{myheadings}
\markboth{\scriptsize  K. Boruah and B. Hazarika}
        {\scriptsize  Some basic properties G-Calculus....}

\maketitle\vspace{-0.4cm}
\section{Introduction}
The  area of non-Newtonian calculus pioneering work carried out  by Grossman and Katz \cite{GrossmanKatz}  which we call as multiplicative calculus. The operations of multiplicative calculus are called as multiplicative derivative and multiplicative integral. We refer to Grossman and Katz \cite{GrossmanKatz}, Stanley \cite{Stanley}, Bashirov et
al. \cite{BashirovMisirh,BashirovKurpinar}, Grossman  \cite{Grossman83} for elements of multiplicative calculus and
its applications. An extension of multiplicative calculus to functions of complex variables is
handled in Bashirov and R\i za \cite{BashirovRiza}, Uzer \cite{Uzer10}, Bashirov et al. \cite{BashirovKurpinar}, \c{C}akmak and Ba\c{s}ar \cite{CakmakBasar}, Cakir \cite{Cakir}, Kadak et al \cite{KadakEfe,kadak2}, Tekin and Ba\c{s}ar\cite{TekinBasar}, T\"{u}rkmen and Ba\c{s}ar \cite{TurkmenBasar}.  Kadak and  \"{O}zl\"{u}k \cite{KADAK3} studied the generalized Runge-Kutta method with
respect to non-Newtonian calculus. 

Geometric calculus is an alternative to the usual calculus of Newton and Leibniz. It provides differentiation and integration tools based on multiplication instead of addition. Every property in Newtonian calculus has an analog in multiplicative calculus. Generally speaking multiplicative calculus is a methodology that allows one to have a different look at problems which can be investigated via calculus. In some cases, for example for growth related problems, the use of multiplicative calculus is advocated instead of a traditional Newtonian one.

\section{$\alpha-$generator and geometric complex field}
A $generator$ is a one-to-one function whose domain is $\mathbb{R}$(the set of real numbers) and whose range is a subset $B\subset \mathbb{R}.$ Each generator generates exactly one arithmetic and each arithmetic is generated by exactly one generator. For example, the identity function generates classical arithmetic, and exponential function generates geometric arithmetic. As a generator, we choose the function $\alpha$ such that whose basic algebraic operations are defined as follows: 
\begin{align*}
&\alpha -addition &x\dot{+}y &=\alpha[\alpha^{-1}(x) + \alpha^{-1}(y)]\\
&\alpha-subtraction &x\dot{-}y&=\alpha[\alpha^{-1}(x) - \alpha^{-1}(y)]\\
&\alpha-multiplication &x\dot{\times}y &=\alpha[\alpha^{-1}(x) \times \alpha^{-1}(y)]\\
&\alpha-division &\dot{x/y}&=\alpha[\alpha^{-1}(x) / \alpha^{-1}(y)]\\
&\alpha-order &x\dot{<}y &\Leftrightarrow \alpha^{-1}(x) < \alpha^{-1}(y).
\end{align*}
for $x, y \in A,$ where $A$ is a domain of the function $\alpha.$

If we choose \textit{$exp$} as an $\alpha-generator$ defined by $\alpha (z)= e^z$ for $z\in \mathbb{C}$ then $\alpha^{-1}(z)=\ln z$ and $\alpha-arithmetic$ turns out to geometric arithmetic.
\begin{align*}
&\alpha -addition &x\oplus y &=\alpha[\alpha^{-1}(x) + \alpha^{-1}(y)]& = e^{(\ln x+\ln y)}& =x.y ~geometric ~addition\\
&\alpha-subtraction &x\ominus y&=\alpha[\alpha^{-1}(x) - \alpha^{-1}(y)]&= e^{(\ln x-\ln y)} &=  x\div y, y\ne 0 ~geometric ~subtraction\\
&\alpha-multiplication &x\odot y &=\alpha[\alpha^{-1}(x) \times\alpha^{-1}(y)]& = e^{(\ln x\times\ln y)} & = ~x^{\ln y}~ geometric ~multiplication\\
&\alpha-division &x\oslash y&=\alpha[\alpha^{-1}(x) / \alpha^{-1}(y)] & = e^{(\ln x\div \ln y)}& = x^{\frac{1}{\ln y}}, y\ne 1 ~ geometric ~division.
\end{align*}
It is obvious that $\ln(x) < \ln(y)$ if $x<y$ for $x, y \in \mathbb{R}^+.$ That is, $x<y \Leftrightarrow \alpha^{-1}(x) < \alpha^{-1}(y)$ So, without loss of generality, we use $x<y$ instead of the geometric order $x\dot{<}y .$

C. T\"{u}rkmen and F. Ba\c{s}ar \cite{TurkmenBasar} defined the sets of geometric integers, geometric real numbers and geometric complex numbers $\mathbb{Z}(G), \mathbb{R}(G)$ and $\mathbb{C}(G),$  respectively, as follows:
\begin{align*}
\mathbb{Z}(G)&=\{ e^{x}: x\in \mathbb{Z}\} = \mathbb{Z}\backslash \{0\}\\
\mathbb{R}(G)&=\{ e^{x}: x\in \mathbb{R}\} = \mathbb{R}^{+}\backslash \{0\}\\
\mathbb{C}(G)&=\{ e^{z}: z\in \mathbb{C}\} = \mathbb{C}\backslash \{0\}.
\end{align*}

\begin{rem}
 $(\mathbb{R}(G), \oplus, \odot)$ is a field with geometric zero $1$ and geometric identity $e,$ since
\begin{itemize}
\item[(a).] $(\mathbb{R}(G), \oplus)$ is a geometric additive Abelian group with geometric zero $1,$
\item[(b).] $(\mathbb{R}(G)\backslash 1, \odot)$ is a geometric multiplicative Abelian group with geometric identity $e,$
\item[(c).] $\odot$ is distributive over $\oplus.$ 
\end{itemize}
\end{rem}

But $(\mathbb{C}(G), \oplus, \odot)$ is not a field, however, geometric binary operation $\odot$ is not associative in $\mathbb{C}(G)$. For, we take
$x = e^{1/4}, y = e^4$ and $z = e^{(1 + i\pi/2)}= ie.$ Then
$(x\odot y)\odot z = e \odot z = z = ie$ but $x\odot(y\odot z) = x\odot e^4 = e.$


Let us define geometric positive real numbers and geometric negative real numbers as follows:
\begin{align*}
\mathbb{R}^+(G)&=\{x\in \mathbb{R}(G) : x >1\}\\
\mathbb{R}^-(G)&=\{x\in \mathbb{R}(G) : x <1\}.
\end{align*}
\subsection{Some useful relations between geometric operations and ordinary arithmetic operations}
For all $x, y\in \mathbb{R}(G)$
\begin{itemize}
\item{ $x\oplus y=xy$}
\item{ $x\ominus y=x/y$}
\item{ $x\odot y=x^{\ln y}=y^{\ln x}$}
\item{ $x\oslash y$ or $\frac{x}{y}G=x^{\frac{1}{\ln y}}, y\neq 1$}
\item{ $x^{2_G}= x \odot x=x^{\ln x}$}
\item{ $x^{p_G}=x^{\ln^{p-1}x}$}
\item{ ${\sqrt{x}}^G=e^{(\ln x)^\frac{1}{2}}$}
\item{ $x^{-1_G}=e^{\frac{1}{\log x}}$}
\item{ $x\odot e=x$ and $x\oplus 1= x$}
\item{ $e^n\odot x=x\oplus x\oplus .....(\text{upto $n$ number of $x$})=x^n$}
\item{
\begin{equation*}
\left|x\right|^G=
\begin{cases}
x, &\text{if $x>1$}\\
1,&\text{if $x=1$}\\
\frac{1}{x},&\text{if $x<1$}
\end{cases}
\end{equation*}}
Thus $\left|x\right|^G\geq 1.$
\item{ ${\sqrt{x^{2_G}}}^G=\left|x\right|^G$}
\item{ $\left|e^y\right|^G=e^{\left|y\right|}$}
\item{ $\left|x\odot y\right|^G=\left|x\right|^G \odot \left|y\right|^G$}
\item{ $\left|x\oplus y\right|^G \leq\left|x\right|^G \oplus \left|y\right|^G$}
\item{ $\left|x\oslash y\right|^G=\left|x\right|^G \oslash \left|y\right|^G$}
\item{ $\left|x\ominus y\right|^G\geq\left|x\right|^G \ominus \left|y\right|^G$}
\item{ $0_G \ominus 1_G\odot\left(x \ominus y\right)=y\ominus x\,, i.e.$ in short $\ominus \left(x \ominus y\right)= y\ominus x.$}
\end{itemize}
Further $e^{-x}=\ominus e^x$ holds for all $x\in \mathbb{Z}^+ .$ Thus the set of all geometric integers turns out to the following:
\[\mathbb{Z}(G)=\{....,e^{-3},e^{-2},e^{-1},e^{0},e^1,e^2,e^3,....\}=\{....,\ominus e^3, \ominus e^2, \ominus e, 1,e,e^2,e^3,....\}.\]
\section{Basic definitions}
\begin{figure}
\centering
\includegraphics[width=0.4\textwidth]{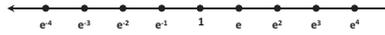}
\caption{Geometric Number Line}
\label{gnl}
\end{figure}

\begin{description}
\item[Geometric Bionomial Formula]
\end{description}
\begin{eqnarray*}
&(i)~~(a\oplus b)^{2_G}&=(a\oplus b)\odot (a \oplus b)\\
                      &&=a\odot a \oplus a\odot b \oplus b\odot a \oplus b\odot b\\
											&&=a^{2_G}\oplus e^2\odot a\odot b\oplus b^{2_G}.\\
&(ii)~~(a\oplus b)^{3_G}&=a^{3_G}\oplus e^3\odot a^{2_G}\odot b\oplus e^3\odot a\odot b^{2_G}\oplus b^{3_G}.\\
\text{In general}&&\\ 
&(iii)~~(a\oplus b)^{n_G}&=a^{n_G}\oplus e^{\binom{n}{1}}\odot a^{(n-1)_G}\odot b\oplus e^{\binom{n}{2}}\odot a^{(n-2)_G}\odot b^{2_G}\oplus....\oplus b^{n_G}\\
                      &&= _G\sum_{r=0}^n e^{\binom{n}{r}}\odot a^{(n-r)_G}\odot b^{r_G}.\\
											\text{Similarly} &&\\
											&~~~~(a\ominus b)^{n_G}&= _G\sum_{r=0}^n \left(\ominus e\right)^{r_G}\odot e^{\binom{n}{r}}\odot a^{(n-r)_G}\odot b^{r_G}.
\end{eqnarray*}
\textit{Note}: $x\oplus x =x^2.$ Also $e^2\odot x= x^{\ln (e^2)}=x^2.$ So, $e^2\odot x = x^2= x\oplus x.$
\begin{description}
\item[Geometric Real Number Line]
\end{description} For $x, y\in \mathbb{R(G)},$ there exist $u, v \in \mathbb{R}$ such that $x=e^u$ and $y=e^v.$ Also consecutive natural numbers are equally spaced by one unit in real number line, but the geometric integers $e, e^2, e^3,...$ are not equally spaced in ordinary sense, e.g. $e^2-e=4.6708$(approx.), $e^3-e^2=12.6965$(approx.). But they are geometrically equidistant as $e^2\ominus e= e^{2-1}=e,e^3\ominus e^2= e^{3-2}=e$ etc. Furthermore, it can be easily verified that $(\mathbb{R(G)}, \oplus, \odot)$ is a complete field with geometric identity $e$ and geometric zero $1.$ So we can consider a new type of geometric real number line as shown in Figure \ref{gnl}.

\begin{figure}
\centering
\includegraphics[width=0.4\textwidth]{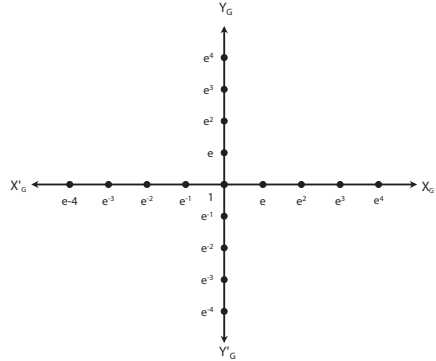}
\caption{Geometric Co-ordinate System}
\label{gcs}
\end{figure}

\begin{figure}
\centering
\includegraphics[width=0.4\textwidth]{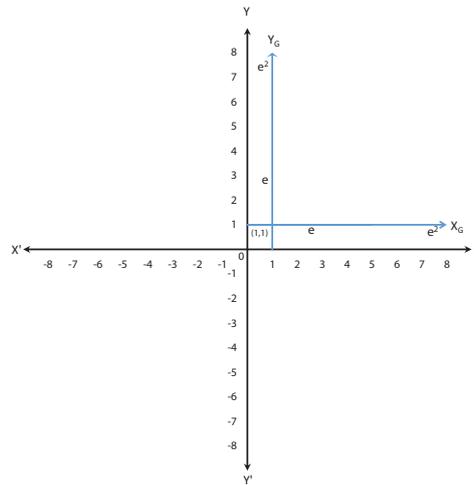}
\caption{Comparison of Geometric axes w.r.t. ordinary axes}
\label{C}
\end{figure}

\begin{description}
\item[Geometric Co-ordinate System]
\end{description}
We consider two mutually perpendicular geometric real number lines which intersect each other at $(1, 1)$ as shown in figure \ref{gcs}. If we compare the geometric axes with respect the ordinary cartesian coordinate system(Figure \ref{C}), the points $e, e^2$ etc. will not be equidistant.

Since consecutive geometric integers are equidistant in geometric sense and $(\mathbb{R}(G), \oplus , \odot)$ is a complete field, so almost all the properties of ordinary cartesian coordinate system will be valid for geometric coordinate system under geometric arithmetic.

\section{G-CALCULUS}
Grossman and Katz\cite{GrossmanKatz} defined the multiplicative differentiation of a function $f(x)$ as 
\begin{equation*}
\frac{d^*f}{dx}=f^*(x)=\lim_{h\rightarrow 0}\left[\frac{f(x+h}{f(x)}\right]^{\frac{1}{h}}
\end{equation*}
We define the G-differentiation of $f(x)$ as
\begin{equation*}
\frac{d^Gf}{dx}=f^G(x)=_G\lim_{h\rightarrow 1}\frac{f(x\oplus h)\ominus f(x)}{h}G \text{~for~} h\in \mathbb{R(G)}.
\end{equation*}
Equivalently
\begin{align*}
\frac{d^Gf}{dx}&=_G\lim_{h\rightarrow 1}\frac{f(x\oplus h)\ominus f(x)}{h}G\\
               &=_G\lim_{h\rightarrow 1}\left[\frac{f(hx)}{f(x)}\right]^{\frac{1}{\ln h}}\\
								&=_G\lim_{u\rightarrow 0}\left[\frac{f(e^u.x)}{f(x)}\right]^{\frac{1}{u}} \text{~where~} h=e^u\in \mathbb{R(G)}\		\end{align*}
The second derivative of $f(x)$ is defined as
\begin{align*}
\frac{d^{2_G}f}{dx^{2_G}}=f^{(2_G)}(x)&= _G\lim_{h\rightarrow 1}\frac{f^G(x\oplus h)\ominus f^G(x)}{h}G
\end{align*} 
Similarly, the $n^{\text{th}}$ derivative is
\begin{align*}
\frac{d^{n_G}f}{dx^{n_G}}=f^{(n_G)}(x)&= _G\lim_{h\rightarrow 1}\frac{f^{\left((n-1)_G\right)}(x\oplus h)\ominus f^{\left((n-1)_G\right)}(x)}{h}G.
\end{align*}
\begin{example}
If $f(x)=x^{n_G},$ then $ f^G(x)= e^n\odot x^{(n-1)_G}$ and $f^{(n_G)}=e^{n!}.$
\end{example}
\begin{proof}
\begin{align*}
f^G(x)&=_G\lim_{h\rightarrow 1}\frac{(x\oplus h)^{n_G}\ominus x^{n_G}}{h}G\\
      &=_G\lim_{h\rightarrow 1}\frac{x^{n_G}\oplus e^{\binom{n}{1}}\odot x^{(n-1)_G}\odot h\oplus e^{\binom{n}{2}}\odot x^{(n-2)_G}\odot h^{2_G}\oplus....\oplus h^{n_G}\ominus x^{n_G}}{h}G\\
			&=_G\lim_{h\rightarrow 1}\frac{e^{\binom{n}{1}}\odot x^{(n-1)_G}\odot h\oplus e^{\binom{n}{2}}\odot x^{(n-2)_G}\odot h^{2_G}\oplus....\oplus h^{n_G}}{h}G\\
			&=_G\lim_{h\rightarrow 1}\left[e^{\binom{n}{1}}\odot x^{(n-1)_G}\oplus e^{\binom{n}{2}}\odot x^{(n-2)_G}\odot h\oplus....\oplus h^{(n-1)_G}\right]\\
			&=e^n\odot x^{(n-1)_G}.\\
\text{Similarly~} f^{(2_G)}(x)&= e^n\odot e^{n-1}\odot x^{(n-2)_G}\\
                            &= e^{n(n-1)}\odot x^{(n-2)_G}.
\end{align*}
Continuing the process, we get ~$f^{(n_G)}(x)= e^{n!}.$
\end{proof}
\begin{example}
Let $f(x)=x.$ Then $\frac{d^Gf}{dx^G}=e.$
\end{example}
\begin{proof} Here
\begin{align*}
 \frac{d^Gf}{dx^G}&= _G\lim_{h\rightarrow 1}\frac{(x\oplus h)\ominus x}{h}G\\
                              &= _G\lim_{h\rightarrow 1}\frac{h}{h}G\\
															&= _G\lim_{h\rightarrow 1}h^{\frac{1}{\ln h}}\\
															&=e \text{~since~}h^{\frac{1}{\ln h}}=e.
\end{align*}
\end{proof}
\section{APPLICATION OF G-CALCULUS IN NUMERICAL ANALYSIS}
\begin{description}
\item[Geometric Factorial]
\end{description}
We defined \cite{KhirodBipan} geometric factorial notation $!_G$ as
\[n!_G=e^n\odot e^{n-1}\odot e^{n-2}\odot \cdots \odot e^2\odot e =e^{n!}.\]
For example,
\begin{align*}
0!_G &=e^{0!}=e^0=1\\
1!_G &=e^{1!}=e=2.71828\\
2!_G &=e^{2!}=e^2=7.38906\\
3!_G &=e^{3!}=e^6=4.03429 \times 10^2\\
4!_G &=e^{4!}=e^{24}=2.64891 \times 10^{10}\\
5!_G &=e^{5!}=e^{120}=1.30418 \times 10^{52}\quad \text{etc.}
\end{align*}
\begin{description}
\item[Generalized Geometric Forward Difference Operator]
\end{description}
Let
\begin{align*}
\Delta_G f(a)   &= f(a\oplus h) \ominus f(a).\\
\Delta^2_G f(a) &= \Delta_G f(a\oplus h) \ominus \Delta_G f(a)\\
                &= \{f(a\oplus e^2\odot h) \ominus f(a\oplus h)\}\ominus \{f(a\oplus h) \ominus f(a)\}\\                
								&=f(a\oplus e^2\odot h) \ominus e^2 \odot f(a\oplus h)\oplus f(a).\\
\Delta^3_G f(a) &= \Delta^2_G f(a\oplus h) \ominus \Delta^2_G f(a)\\
                &=\{f(a\oplus e^3\odot h) \ominus e^2 \odot f(a\oplus e^2\odot h)\oplus f(a \oplus h)\}\\
								&\qquad \ominus \{f(a\oplus e^2\odot h) \ominus e^2 \odot f(a\oplus h)\oplus f(a)\}\\
								&=f(a\oplus e^3\odot h) \ominus e^3 \odot f(a\oplus e^2\odot h)\oplus e^3\odot f(a \oplus h) \ominus f(a).
\end{align*}
Thus, $n^{\text{th}}$ forward difference is
\[\Delta^n_G f(a)= _G\sum^n_{k=0} (\ominus e )^{{k}_G}\odot e^{\binom{n}{k}}\odot f(a\oplus e^{n-k}\odot h), \text{with}\, (\ominus e)^{0_G}=e.\]\\

\begin{description}
\item[Generalized Geometric Backward Difference Operator]
\end{description}
 Let
\begin{align*}
\nabla_G f(a)   &=f(a) \ominus f(a\ominus h).\\
\nabla^2_G f(a) &= \nabla_G f(a) \ominus \nabla_G f(a\ominus h)\\
                &= \{f(a)\ominus f(a \ominus h)\} \ominus \{f(a\ominus h) \ominus f(a\ominus e^2\odot h)\}\\
								&= f(a)\ominus  e^2 \odot f(a \ominus h)\oplus f(a \ominus e^2 \odot h).\\
\nabla^3_G f(a) &= \nabla^2_G f(a) \ominus \nabla^2_G f(a-h)\\
                &=\{f(a)\ominus  e^2 \odot f(a \ominus h)\oplus f(a \ominus e^2 \odot h)\}\\
								&\qquad \ominus \{f(a\ominus h)\ominus  e^2 \odot f(a \ominus e^2 \odot h)\oplus f(a \ominus e^3 \odot h)\}\\
								&=f(a) \ominus e^3 \odot f(a \ominus h)\oplus e^3\odot f(a \ominus e^2 \odot h) \ominus f(a\ominus e^3 \odot h).
\end{align*}
Thus, $n^{\text{th}}$ geometric backward difference is
\[\nabla^n_G f(a)= _G\sum^n_{k=0} (\ominus e )^{{k}_G}\odot e^{\binom{n}{k}}\odot f(a\ominus e^k\odot h).\]
\begin{description}
\item[Divided difference]
\end{description}
 Let $f(x_0), f(x_1),...,f(x_n)$ be the entries corresponding to the arguments $x_0, x_1,...,x_n$ where the intervals $x_1\ominus x_0, x_2\ominus x_1,..., x_n\ominus x_{n-1}$ not necessarily equal, i.e. values of the argument are not geometrically equally spaced. Then the first divided difference of $f(x)$ for the arguments $x_0, x_1$ is defined as $\frac{f(x_1)\ominus f(x_0)}{x_1\ominus x_0}G$ or $\frac{f(x_0)\ominus f(x_1)}{x_0\ominus x_1}G$ and is denoted by $\underset{x_1}{\blacktriangle_G}f(x_0)$ or $f(x_0, x_1).$ i.e.
\begin{align*}
f(x_0, x_1)&= \frac{f(x_1)\ominus f(x_0)}{x_1\ominus x_0}G=\left[\frac{f(x_1)}{f(x_0)}\right]^{\frac{1}{\ln{(\frac{x_1}{x_0})}}}=\underset{x_1}{\blacktriangle_G}f(x_0).\text{~Similarly}\\
f(x_1, x_2)&= \frac{f(x_2)\ominus f(x_1)}{x_2\ominus x_1}G=\left[\frac{f(x_2)}{f(x_1)}\right]^{\frac{1}{\ln{(\frac{x_2}{x_1})}}}=\underset{x_2}{\blacktriangle_G}f(x_1)\\
f(x_2, x_3)&= \frac{f(x_3)\ominus f(x_2)}{x_3\ominus x_2}G=\left[\frac{f(x_3)}{f(x_2)}\right]^{\frac{1}{\ln{(\frac{x_3}{x_2})}}}=\underset{x_3}{\blacktriangle_G}f(x_2)\\
\hdots \hspace{1cm}&\hdots \hspace{1cm}\hdots \hspace{1cm}\hdots \hspace{1cm}\hdots\\
\hdots \hspace{1cm}&\hdots \hspace{1cm}\hdots \hspace{1cm}\hdots \hspace{1cm}\hdots\\
f(x_{n-1}, x_n)&= \frac{f(x_n)\ominus f(x_{n-1})}{x_n\ominus x_{n-1}}G=\left[\frac{f(x_n)}{f(x_{n-1})}\right]^{\frac{1}{\ln{(\frac{x_n}{x_{n-1}})}}}=\underset{x_n}{\blacktriangle_G}f(x_{n-1}).
\end{align*}
The second geometric divided difference of $f(x)$ for the three arguments $x_0, x_1$ and $x_2$ is defined as
\begin{align*}
f(x_0, x_1,x_2)&= \frac{f(x_1, x_2)\ominus f(x_0, x_1)}{x_2\ominus x_0}G=\left[\frac{f(x_1, x_2)}{f(x_0, x_1)}\right]^{\frac{1}{\ln{(\frac{x_2}{x_1})}}}=\underset{x_1, x_2}{\blacktriangle^2_G}f(x_0).
\end{align*}
The $n^{\text{th}}$ divided difference is given by
\begin{align*}
f(x_0, x_1,...,x_n)&= \frac{f(x_1, x_2,...,x_n)\ominus f(x_0, x_1,...,x_{n-1})}{x_n\ominus x_0}G\\
                   &=\left[\frac{f(x_1, x_2,...,x_n)}{f(x_0, x_1,...,x_{n-1})}\right]^{\frac{1}{\ln{(\frac{x_n}{x_0})}}}\\
									 &=\underset{x_1, x_2,...,x_n}{\blacktriangle^n_G}f(x_0).
\end{align*}
For convenience, we'll write ${\blacktriangle^n_G}f(x_0)$ instead of $\underset{x_1, x_2,...,x_n}{\blacktriangle^n_G}f(x_0).$
\begin{rem}\label{rem1}
If two of the arguments coincide, the divided difference can be given by taking limit as:
\begin{align*}
f(x_0, x_0)=\lim_{\epsilon \to 1}f(x_0, x_0\oplus \epsilon)&=\lim_{\epsilon \to 1}\frac{f(x_0\oplus \epsilon)\ominus f(x_0)}{\epsilon}G\\
&=f^G(x_0),\text{~which is the geometric derivative of $f(x)$ at $x_0$}.
\end{align*}
Similarly, $f(\underbrace{x_0,x_0,x_0,...,x_0}_{(m+1)\text~{arguments}})=\frac{e}{e^{m!}}G\odot f^{(m_G)}(x_0).$
\end{rem} 
\begin{rem}\label{rem2}
The $n^{\text{th}}$ divided differences of a geometric polynomial of degree $n$ are constant. 
\end{rem}
\begin{proof}
First we consider a function $f(x)=x^{n_G},$ i.e. in ordinary sense $f(x)=x^{\ln^{n-1}x}.$Then the first divided differences of this function are given by
\begin{align*}
f(x_r,x_{r+1})&=\frac{f(x_{r+1})\ominus f(x_r)}{x_{r+1}\ominus x_r}G=\frac{x^{n_G}_{r+1}\ominus x^{n_G}_r}{x_{r+1}\ominus x_r}\\
              &=x^{(n-1)_G}_{r+1}\oplus x^{(n-2)_G}_{r+1}\odot x_r\oplus x^{(n-3)_G}_{r+1}\odot x^{2_G}_r\oplus\hdots \oplus x_{r+1}\odot x^{(n-2)_G}_r\oplus x^{(n-1)_G}_r,
\end{align*}
which is a homogeneous expression of degree $(n-1)$ in $x_r$ and $x_{r+1}.$ The second divided differences are given by
\begin{gather}
\begin{aligned}
f(x_r,x_{r+1},x_{r+2})&=\frac{f(x_{r+1}, x_{r+2})\ominus f(x_r, x_{r+1})}{x_{r+2}\ominus x_r}G\\
              &=\frac{x^{(n-1)_G}_{r+2}\oplus x^{(n-2)_G}_{r+2}\odot x_{r+1}\oplus x^{(n-3)_G}_{r+2}\odot x^{2_G}_{r+1}\oplus\hdots \oplus x^{(n-1)_G}_{r+1}}{x_{r+2}\ominus x_r}G\\
							&\ominus \frac{x^{(n-1)_G}_{r+1}\oplus x^{(n-2)_G}_{r+1}\odot x_r\oplus x^{(n-3)_G}_{r+1}\odot x^{2_G}_r\oplus\hdots \oplus x^{(n-1)_G}_r}{x_{r+2}\ominus x_r}G\\
							&= \frac{x^{(n-1)_G}_{r+2}\ominus x^{(n-1)_G}_r}{x_{r+2}\ominus x_r}G\oplus x_{r+1}\odot\frac{x^{(n-2)_G}_{r+2}\ominus x^{(n-2)_G}_r}{x_{r+2}\ominus x_r}G\oplus \\
							&\qquad\hdots \oplus x^{(n-2)_G}_{r+1}\odot \frac{x_{r+2}\ominus x_r}{x_{r+2}\ominus x_r}G\\
							&=\left(x^{(n-2)_G}_{r+2}\oplus \hdots \oplus x^{(n-2)_G}_r\right)\oplus x_{r+1}\odot \left(x^{(n-3)_G}_{r+2}\oplus \hdots \oplus x^{(n-3)_G}_r\right)\oplus \\
							&\qquad\hdots \oplus x^{(n-2)_G}_{r+1},
\end{aligned}
\end{gather}
which is a homogeneous expression of degree $n-2$ in $x_r, x_{r+1}$ and $x_{r+2}.$

By induction it can be shown that the $n^{\text{th}}$ divided difference of $f(x)=x^{n_G}$ is an expression of degree zero, i.e. a constant, and therefore independent of the values $x_r, x_{r+1},...,x_{r+n}.$ Since the $n^{\text{th}}$ geometric divided difference of $x^{n_G}$ are constant, therefore the geometric divided differences of $x^{n_G}$ of order less than $n$ will all be zero.

 If $f(x)=a\odot x^{n_G},$ where $a\in \mathbb{R}(G)$ is a constant, then the $n^{\text{th}}$ geometric divided difference of $f(x)=a\odot (n^{\text{th}} \text{~geometric divided difference of~} x^{n_G}),$ which is a constant.

Therefore if $f(x)=a_0\odot x^{n_G}\oplus a_1\odot x^{(n-1)_G}\oplus\hdots \oplus a_{n-1}\odot x\oplus a_n$ be a geometric polynomial of degree $n,$ so that $a_0, a_1,...,a_n\in \mathbb{R}(G),$ then the $n^{\text{th}}$ geometric divided difference of all terms will be vanished except that $a_0\odot x^n$ turns to constant. Hence $n^{\text{th}}$ geometric divided difference of the whole polynomial will be constant.
\end{proof}
In our paper \cite{KhirodBipan}, we have derived Geometric Newton-Gregory forward interpolation and backward interpolation formulae for geometrically equidistant values of the argument. Here, we derive divided difference interpolation formula for unequal values of the argument.

\begin{description}
\item[Divided difference Interpolation formula]
\end{description}
 Let $f(x_0),f(x_1),f(x_2),...,f(x_n)$ be the values of $f(x)$ corresponding to the arguments $x_0,x_1,x_2,...,x_n$ not necessarily geometrically equally spaced. From the definition of divided differences
\begin{align*}
f(x, x_0)&= \frac{f(x)\ominus f(x_0)}{x\ominus x_0}G\\
(x\ominus x_0)\odot f(x)&= f(x)\ominus f(x_0)\\
\implies f(x)&= f(x_0)\oplus (x\ominus x_0)\odot f(x)\tag{1}\\
\text{Also~~}~ f(x, x_0, x_1)&=\frac{f(x, x_0)\ominus f(x_0, x_1)}{x\ominus x_1}G\\
\implies f(x, x_0) &= f(x_0, x_1)\oplus (x\ominus x_1)\odot f(x, x_0, x_1)\\ 
\therefore \text{~from~} (1), f(x)&= f(x_0)\oplus (x\ominus x_0)\odot \left\{f(x_0,x_1)\oplus (x-x_1)\odot f(x, x_0, x_1)\right\}\\
&=f(x_0)\oplus (x\ominus x_0)\odot f(x_0,x_1)\oplus (x\ominus x_0)(x-x_1)\odot f(x, x_0, x_1)\\
\text{Again~~}~ f(x, x_0, x_1)&= f(x_0, x_1, x_2)\oplus (x-x_2)\odot f(x, x_0, x_1, x_2)\\
\text{Similarly~~}~f(x, x_0, x_1, x_2)&= f(x_0, x_1, x_2)\oplus (x-x_3)\odot f(x, x_0, x_1, x_2, x_3)\\
\hdots \hspace{1cm}&\hdots \hspace{1cm}\hdots \hspace{1cm}\hdots \hspace{1cm}\hdots\\
\hdots \hspace{1cm}&\hdots \hspace{1cm}\hdots \hspace{1cm}\hdots \hspace{1cm}\hdots\\
f(x, x_0, x_1,...,x_{n-1})&= f(x_0, x_1,... x_n)\oplus \hdots \oplus (x-x_n)\odot f(x, x_0, x_1,..., x_n).
\end{align*}
Substituting successively we get
\begin{multline}\label{eqn5}
f(x)=f(x_0)\oplus (x\ominus x_0)\odot f(x_0,x_1)\oplus (x\ominus x_0)\odot (x\ominus x_1)\odot f(x_0,x_1, x_2)\oplus\\ ... \oplus (x\ominus x_0)\odot (x\ominus x_1)\odot ...\odot (x\ominus x_{n-1}) \odot f(x_0,x_1,... ,x_n)\oplus R_n,
\end{multline}
where the the remainder $R_n$ is given by
\begin{align*}
R_n=(x\ominus x_0)\odot (x\ominus x_1)\odot ...\odot(x\ominus x_n)\odot f(x,x_0,x_1,... ,x_n).
\end{align*}
If $f(x)$ is a geometric polynomial of degree n, then $f(x,x_0,x_1,... ,x_n)=1,$ so equation (\ref{eqn5}) becomes
\begin{multline}\label{eqn6}
f(x)=f(x_0)\oplus (x\ominus x_0)\odot f(x_0,x_1)\oplus (x\ominus x_0)\odot (x\ominus x_1)\odot f(x_0,x_1, x_2)\oplus\\ ... \oplus (x\ominus x_0)\odot (x\ominus x_1)\odot ...\odot (x\ominus x_{n-1})\odot f(x_0,x_1,... ,x_n).
\end{multline}
This is the divided difference interpolation formula for geometrically unequal intervals.
\begin{description}
\item[Relation between divided differences and geometric forward interpolation ]
\end{description}
Let the values of the argument $x_0, x_1, x_2,...,x_n$ be equally spaced, i.e., $x_1=x_0\oplus h, x_2=x_0\oplus e^2\odot h, x_3=x_0\oplus e^3\odot h,...,x_n=x_0\oplus e^{n-1}\odot h$ and let $\frac{x\ominus x_0}{h}G=u.$ Then 
\begin{align*}
f(x_0, x_1)&=\frac{f(x_1)\ominus f(x_0)}{x_1 \ominus x_0}G\\
           &=\frac{f(x_0\oplus h)\ominus f(x_0)}{x_0\oplus h \ominus x_0}G\\
					&= \frac{\Delta_Gf(x_0)}{h}G\\
f(x_0, x_1, x_2)&=\frac{f(x_1, x_2)\ominus f(x_1, x_0)}{x_2 \ominus x_0}G\\
                &= \frac{\frac{\Delta_Gf(x_1)}{h}G\ominus \frac{\Delta_Gf(x_0)}{h}G}{e^2\odot h}G\\
								&=\frac{\Delta^2_Gf(x_0)}{e^2\odot h\odot h}G\\
								&=\frac{\Delta^2_Gf(x_0)}{2!_G\odot h^{2_G}}G.\\
\text{Similarly,}~~~f(x_0, x_1, x_2, x_3)&= \frac{\Delta^3_Gf(x_0)}{3!_G\odot h^{3_G}}G.\\
\hdots \hspace{1cm}&\hdots \hspace{1cm}\hdots \hspace{1cm}\hdots \hspace{1cm}\hdots\\
\hdots \hspace{1cm}&\hdots \hspace{1cm}\hdots \hspace{1cm}\hdots \hspace{1cm}\hdots\\
f(x_0, x_1, x_2,...,x_n)&= \frac{\Delta^n_Gf(x_0)}{n!_G\odot h^{n_G}}G.
\end{align*}
Substituting these values of divided differences in (\ref{eqn6}) we get
\begin{gather}
\begin{aligned}\label{eqn7}
f(x)=& f(x_0)\oplus u\odot \Delta_G f(x_0) \oplus \frac{u\odot(u\ominus e)}{2!_G}G\odot \Delta^2_G f(x_0)\\
      &\oplus \frac{u\odot(u\ominus e)\odot (u \ominus e^2)}{3!_G}G \odot \Delta^3_G f(x_0)\oplus\cdots\\
			&\oplus \frac{u\odot (u \ominus e)\odot(u \ominus e^2)\odot \cdots \odot (u \ominus e^{n-1})}{n!_G}G \odot \Delta^n_G f(x_0).
\end{aligned}
\end{gather}
This is the Newton-Gregory formula for geometric forward interpolation about which has been discussed in \cite{KhirodBipan}.

\begin{rem}
The geometric divided differences are symmetrical in all their arguments, i.e. the value of any difference is independent of the order of the arguments.
\end{rem}
\begin{proof}
Let $f(x_0),f(x_1),f(x_2),...,f(x_n)$ be the values of $f(x)$ corresponding to the arguments $x_0,x_1,x_2,...,x_n.$ We have the first geometric divided difference (GDD) is
\begin{align*}
f(x_0, x_1)&= \frac{f(x_1)\ominus f(x_0)}{x_1\ominus x_0}G= \frac{f(x_0)\ominus f(x_1)}{x_0\ominus x_1}G=f(x_1, x_0)\\
           &=\frac{f(x_0)}{x_0\ominus x_1}G\oplus \frac{f(x_1)}{x_1\ominus x_0}G=_G\sum\frac{f(x_0)}{x_0\ominus x_1}G,
\end{align*}
showing that $f(x_0, x_1)$ is symmetrical in $x_0,x_1.$ Again, the second GDD is
\begin{align*}
f(x_0, x_1,x_3)&= \frac{f(x_1,x_2)\ominus f(x_0,x_1)}{x_2\ominus x_0}G\\
               &=\frac{e}{x_2\ominus x_0}G\odot\left[\left\{\frac{f(x_1)}{x_1\ominus x_2}G\oplus \frac{f(x_2)}{x_2\ominus x_1}G\right\}\ominus\left\{\frac{f(x_0)}{x_0\ominus x_1}G\oplus \frac{f(x_1)}{x_1\ominus x_0}G\right\}\right]\\
							&=\frac{f(x_0)}{(x_0\ominus x_1)\odot (x_0\ominus x_2)}G\oplus \frac{f(x_1)}{(x_1\ominus x_0)\odot (x_1\ominus x_2)}G\\
							&\qquad\oplus \frac{f(x_2)}{(x_2\ominus x_0)\odot (x_2\ominus x_1)}G\\
							&=_G\sum \frac{f(x_0)}{(x_0\ominus x_1)\odot (x_0\ominus x_2)}G,\text{~showing that second GDD is symmetrical}.
\end{align*}

Let us assume that $(n-1)^{\text{th}}$ GDD is symmetrical, i.e.
\begin{equation*}
\begin{split}
f(x_0, x_1,...,x_{n-1})&= \frac{f(x_0)}{(x_0\ominus x_1)\odot ...\odot (x_0\ominus x_{n-1})}G\\
              &\qquad\oplus\frac{f(x_1)}{(x_1\ominus x_0)\odot (x_1\ominus x_2)\odot ...\odot (x_1\ominus x_{n-1})}G\oplus ...\\
							&\qquad\oplus\frac{f(x_{n-1})}{(x_{n-1}\ominus x_0)\odot (x_{n-1}\ominus x_1)\odot ...\odot (x_{n-1}\ominus x_{n-2})}G\\
							&= _G\sum\frac{f(x_0)}{(x_0\ominus x_1)\odot ...\odot (x_0\ominus x_{n-1})}G
\end{split}
\end{equation*}
Then
\begin{equation}\label{eqn8}
\begin{split}
f(x_0, x_1,...,x_n)&= \frac{f(x_0,...,x_{n-1})\ominus f(x_1,..,x_n)}{x_0\ominus x_n}G\\
                   &=\frac{e}{x_0\ominus x_n}G\odot\left\{_G\sum\frac{f(x_0)}{(x_0\ominus x_1)\odot ...\odot (x_0\ominus x_{n-1})}G\right\}\\
									&\qquad \ominus \frac{e}{x_0\ominus x_n}G\odot\left\{_G\sum\frac{f(x_1)}{(x_1\ominus x_0)\odot(x_1\ominus x_2)\odot ...\odot (x_1\ominus x_n)}G\right\}\\
                   &= \frac{f(x_0)}{(x_0\ominus x_1)\odot ...\odot (x_0\ominus x_n)}G\\
              &\qquad \oplus\frac{f(x_1)}{(x_1\ominus x_0)\odot (x_1\ominus x_2)\odot ...\odot (x_1\ominus x_n)}G\oplus ...\\
							&\qquad\oplus\frac{f(x_n)}{(x_n\ominus x_0)\odot (x_n\ominus x_1)\odot ...\odot (x_n\ominus x_{n-1})}G\\
							    & = _G\sum\frac{f(x_0)}{(x_0\ominus x_1)\odot ...\odot (x_0\ominus x_n)}G
\end{split}
\end{equation}
\end{proof}
\begin{description}
\item[Lagrange's Geometric Interpolation Formula for Unequal Intervals]
\end{description}
Let $y=f(x)$ be a geometric polynomial of degree $n$ which takes the values $y_0=f(x_0), y_1=f(x_1),...,y_n=f(x_n)$ as $x$ takes the values $x_0, x_1,...,x_n,$ respectively. Then the $(n+1)^{\text{th}}$ GDD of this polynomial are zero. Hence
 \begin{equation}\label{eqn9}
f(x,x_0, x_1,x_2,...,x_n)=0
\end{equation}
Again from equation (\ref{eqn8})
\begin{equation}
\begin{split}
f(x,x_0, x_1,x_2,...,x_n)&=\frac{y}{(x\ominus x_0)\odot (x\ominus x_1)\odot (x\ominus x_2)\odot ...\odot (x\ominus x_n)}G\\
                          &\qquad\oplus \frac{y_0}{(x_0\ominus x)\odot (x_0\ominus x_1)\odot (x_0\ominus x_2)\odot ...\odot (x_0\ominus x_n)}G\\
													&\qquad\oplus \frac{y_1}{(x_1\ominus x)\odot (x_1\ominus x_0)\odot (x_1\ominus x_2)\odot ...\odot (x_1\ominus x_n)}G\oplus ...\\
													&\qquad\oplus \frac{y_n}{(x_n\ominus x)\odot (x_n\ominus x_0)\odot (x_n\ominus x_0)\odot ...\odot (x_n\ominus x_{n-1})}G
\end{split}
\end{equation}
Using (\ref{eqn9}), we get
\begin{equation}
\begin{split}
&\frac{y}{(x\ominus x_0)\odot (x\ominus x_1)\odot (x\ominus x_2)\odot ...\odot (x\ominus x_n)}G\\
                          &\quad\oplus \frac{y_0}{(x_0\ominus x)\odot (x_0\ominus x_1)\odot (x_0\ominus x_2)\odot ...\odot (x_0\ominus x_n)}G\\
													&\qquad\oplus \frac{y_1}{(x_1\ominus x)\odot (x_1\ominus x_0)\odot (x_1\ominus x_2)\odot ...\odot (x_1\ominus x_n)}G\oplus ...\\
													&\qquad \quad \oplus \frac{y_n}{(x_n\ominus x)\odot (x_n\ominus x_0)\odot (x_n\ominus x_0)\odot ...\odot (x_n\ominus x_{n-1})}G=0
\end{split}
\end{equation}
Transposing all terms except the first term, to the right hand side, we get
\begin{equation}
\begin{split}
&\frac{y}{(x\ominus x_0)\odot (x\ominus x_1)\odot (x\ominus x_2)\odot ...\odot (x\ominus x_n)}G\\
                          &=\oplus \frac{y_0}{(x\ominus x_0)\odot (x_0\ominus x_1)\odot (x_0\ominus x_2)\odot ...\odot (x_0\ominus x_n)}G\\
													&\qquad\oplus \frac{y_1}{(x\ominus x_1)\odot (x_1\ominus x_0)\odot (x_1\ominus x_2)\odot ...\odot (x_1\ominus x_n)}G\oplus ...\\
													&\qquad \quad \oplus \frac{y_n}{(x\ominus x_n)\odot (x_n\ominus x_0)\odot (x_n\ominus x_0)\odot ...\odot (x_n\ominus x_{n-1})}G
\end{split}
\end{equation}
Multiplying both sides by $(x\ominus x_0)\odot (x\ominus x_1)\odot (x\ominus x_2)\odot ...\odot (x\ominus x_n),$ we get
\begin{equation}
\begin{split}
y   &=\frac{(x\ominus x_1)\odot (x\ominus x_2)\odot ...\odot (x\ominus x_n)}{(x_0\ominus x_1)\odot (x_0\ominus x_2)\odot ...\odot (x_0\ominus x_n)}G\odot y_0\\
    &\qquad\oplus \frac{(x\ominus x_0)\odot (x\ominus x_2)\odot ...\odot (x\ominus x_n)}{(x_1\ominus x_0)\odot (x_1\ominus x_2)\odot ...\odot (x_1\ominus x_n)}G\odot y_1\\
		&\qquad \quad \oplus \frac{(x\ominus x_0)\odot (x\ominus x_1)\odot (x\ominus x_3)\odot ...\odot (x\ominus x_n)}{(x_2\ominus x_0)\odot (x_2\ominus x_1)\odot (x_2\ominus x_3)\odot ...\odot (x_1\ominus x_n)}G\odot y_2\oplus ...\\
    &\qquad \qquad \oplus \frac{(x\ominus x_0)\odot (x\ominus x_1)\odot (x\ominus x_2)\odot ...\odot (x\ominus x_{n-1})}{(x_n\ominus x_0)\odot (x_n\ominus x_1)\odot (x_n\ominus x_2)\odot ...\odot (x_n\ominus x_{n-1})}G\odot y_n
\end{split}
\end{equation}
\begin{example}\label{example1} Given, $f(x)=f(e^t)=\sin(e^t).$ From the following table, find $\sin(e^{0.14})$ using geometric divided difference formula.\\[2ex]
\begin{center}
\centering 
\begin{tabular}{|c| c| c| c|c|} 
\hline
$x$  & $e^{0.12}$ &$e^{0.15}$  & $e^{0.19}$& $e^{0.21}$\\ [0.5ex]
\hline 
$f(x)$&$0.903341$& $0.917534$&$0.935351$&$0.943712$\\[1ex] 
\hline 
\end{tabular}\\[3ex]
\end{center}
\textbf{Solution:}
The geometric divided difference table for given geometrically unequal data is as follows:\\[2ex]
\begin{center}
\centering 
\begin{tabular}{|c| c| c| c| c|} 
\hline
$x$(in radian)       & $f(x)$ & $\blacktriangle_G f(x)$  & $\blacktriangle^2_G f(x)$& $\blacktriangle^3_G f(x)$\\ [1.5ex]
\hline 
$e^{0.12}$& 0.903341 &                  &                  & \\
          &          &1.681421          &                  & \\
$e^{0.15}$& 0.917534 &                  & 0.574158         & \\
          &          &1.617367          &                  & 0.623266 \\
$e^{0.19}$& 0.935351 &                  & 0.55024          & \\
          &          &1.560421          &                  & \\
$e^{0.21}$& 0.943712 &                  &                  & \\
\hline 
\end{tabular}\\[1ex]
\end{center}
It is to be noted that here in the table,
\begin{align*}
\blacktriangle_G f(x_0)&=\left[\frac{0.917534}{0.903341}\right]^{\frac{1}{\ln\left({\frac{e^{0.15}}{e^{0.12}}}\right)}}=\left(1.0157116748\right)^{\frac{1}{0.03}}=1.681421\\
\blacktriangle_G f(x_1)&=\left[\frac{0.935351}{0.917534}\right]^{\frac{1}{\ln\left({\frac{e^{0.19}}{e^{0.15}}}\right)}}=1.617367\text{~and so on.}\\
\text{Similarly~} \blacktriangle^2_G f(x_0)&=\left[\frac{1.617367 }{1.681421 }\right]^{\frac{1}{\ln\left({\frac{e^{0.19}}{e^{0.12}}}\right)}}=0.574158 \text{~etc.}
\end{align*}
Now, using the geometric divided difference formula, we get
\begin{align*}
f(x)&=f(x_0)\oplus (x\ominus x_0)\odot \blacktriangle_G f(x_0) \oplus (x\ominus x_0)\odot (x \ominus x_1)\odot \blacktriangle^2_G f(x_0)\\
       &\quad\oplus (x \ominus x_0)\odot (x \ominus x_1)\odot (x \ominus x_2)\odot \blacktriangle^3_G f(x_0)\oplus ...\\
f(e^{0.14}) &= 0.903341\oplus (e^{0.14}\ominus e^{0.12})\odot 1.681421 \\
       &\quad \oplus (e^{0.14}\ominus e^{0.12})\odot(e^{0.14}\ominus e^{0.15})\odot 0.574158\\
			 & \quad \oplus (e^{0.14}\ominus e^{0.12})\odot (e^{0.14}\ominus e^{0.12})\odot (e^{0.15}\ominus e^{0.19})\odot 0.623266\\
			&= 0.903341\oplus e^{0.14-0.12}\odot 1.681421 \oplus (e^{0.14 -0.12}\odot e^{0.14-0.15})\odot 0.574158\\
			 & \quad \oplus (e^{0.14-0.12}\odot e^{0.14-0.12}\odot e^{0.15-0.19})\odot 0.623266\\
			&=0.903341\oplus e^{0.02}\odot 1.681421 \oplus (e^{0.02}\odot e^{-0.01})\odot 0.574158\\
			 & \quad \oplus (e^{0.02}\odot e^{-0.01}\odot e^{-0.05})\odot 0.623266\\
			&=0.903341\oplus e^{0.02}\odot 1.681421 \oplus e^{0.02(-0.01)}\odot 0.574158 \oplus e^{0.02(-0.01)(-0.05)}\odot 0.623266\\
			&=0.903341\oplus e^{0.02}\odot 1.681421 \oplus e^{-0.0002}\odot 0.574158 \oplus e^{0.00001}\odot 0.623266\\
		  &=0.903341\oplus (1.681421)^{0.02} \oplus (0.574158)^{\frac{1}{0.0002}} \oplus (0.623266)^{0.00001}\\
			&=0.903341\oplus 1.010447 \oplus 1.000111 \oplus 0.999995\\
			&=0.903341\times 1.010447 \times 1.000111 \times 0.999995\\
			&=0.912875
\end{align*}
Hence $\sin(e^{0.14})=0.912875$ which is accurate upto the last place of decimal (i.e. sixth place). Thus, geometric divided difference interpolation formula gives values of transcendental functions at given points upto desired degree of accuracy. 
\end{example}
\begin{example}
Solve the problem explained in Example \ref{example1} by Lagrange's geometric interpolation formula.
\end{example}
\textbf{Solution:} Given
\begin{center}
\centering 
\begin{tabular}{|c| c| c| c|c|} 
\hline
$x$  & $e^{0.12}$ &$e^{0.15}$  & $e^{0.19}$& $e^{0.21}$\\ [0.5ex]
\hline 
$f(x)$&$0.903341$& $0.917534$&$0.935351$&$0.943712$\\[1ex] 
\hline 
\end{tabular}\\[3ex]
\end{center}
Here $x= e^{0.14}, x_0=e^{0.12}, x_1=e^{0.15}, x_2=e^{0.19}$ and $x_3=e^{0.21}.$
Now, Lagrange's interpolation formula gives
\begin{equation*}
\begin{split}
f(x) =&\frac{(x\ominus x_1)\odot (x\ominus x_2)\odot(x\ominus x_3)}{(x_0\ominus x_1)\odot (x_0\ominus x_2)\odot(x_0\ominus x_3)}G\odot f(x_0)\\
    &\oplus \frac{(x\ominus x_0)\odot (x\ominus x_2)\odot(x\ominus x_3)}{(x_1\ominus x_0)\odot (x_1\ominus x_2)\odot(x_1\ominus x_3)}G\odot f(x_1)\\
		&\oplus \frac{(x\ominus x_0)\odot (x\ominus x_1)\odot(x\ominus x_3)}{(x_2\ominus x_0)\odot (x_2\ominus x_1)\odot (x_2\ominus x_3)}G\odot f(x_2)\\
    &\oplus \frac{(x\ominus x_0)\odot (x\ominus x_1)\odot (x\ominus x_2)}{(x_3\ominus x_0)\odot (x_3\ominus x_1)\odot (x_3\ominus x_2)}G\odot f(x_3)
\end{split}
\end{equation*}
Putting the values we get
\begin{equation*}
\begin{split}
f(e^{0.14}) =&\frac{(e^{0.14}\ominus e^{0.15})\odot (e^{0.14}\ominus e^{0.19})\odot(e^{0.14}\ominus e^{0.21})}{(e^{0.12}\ominus e^{0.15})\odot (e^{0.12}\ominus e^{0.19})\odot(e^{0.12}\ominus e^{0.21})}G\odot f(e^{0.12})\\
    &\oplus \frac{(e^{0.14}\ominus e^{0.12})\odot (e^{0.14}\ominus e^{0.19})\odot(e^{0.14}\ominus e^{0.21})}{(e^{0.15}\ominus e^{0.12})\odot (e^{0.15}\ominus e^{0.19})\odot(e^{0.15}\ominus e^{0.21})}G\odot f(e^{0.15})\\
		&\oplus \frac{(e^{0.14}\ominus e^{0.12})\odot (e^{0.14}\ominus e^{0.15})\odot(e^{0.14}\ominus e^{0.21})}{(e^{0.19}\ominus e^{0.12})\odot (e^{0.19}\ominus e^{0.15})\odot (e^{0.19}\ominus e^{0.21})}G\odot f(e^{0.19})\\
    &\oplus \frac{(e^{0.14}\ominus e^{0.12})\odot (e^{0.14}\ominus e^{0.15})\odot (e^{0.14}\ominus e^{0.19})}{(e^{0.21}\ominus e^{0.12})\odot (e^{0.21}\ominus e^{0.15})\odot (e^{0.21}\ominus e^{0.19})}G\odot f(e^{0.21})\\
		=&e^{\frac{(-0.01)(-0.05)(-0.07)}{(-0.03)(-0.07)(-0.09)}}\odot 0.903341 \oplus e^{\frac{(0.02)(-0.05)(-0.07)}{(0.03)(-0.04)(-0.06)}}\odot 0.917534\\
		&\oplus e^{\frac{(0.02)(-0.01)(-0.07)}{(0.07)(0.04)(-0.02)}}\odot 0.935351 \oplus e^{\frac{(0.02)(-0.01)(-0.05)}{(0.09)(0.06)(0.02)}} \odot 0.943712\\
		=& e^{\frac{5}{27}}\odot 0.903341  \oplus e^{\frac{35}{36}}\odot 0.917534\oplus e^{\frac{-1}{4}}\odot 0.935351 \oplus e^{\frac{5}{54}}\odot 0.943712\\
		=& (0.903341)^{\frac{5}{27}}  \oplus (0.917534)^{\frac{35}{36}}\oplus (0.935351)^{\frac{-1}{4}} \oplus (0.943712)^{\frac{5}{54}}\\
		=& 0.981351 \oplus 0.91973 \oplus 1.016849 \oplus 0.99465\\
		=& 0.981351 \times 0.91973 \times 1.016849 \times 0.99465\\
		=& 0.912875
\end{split}
\end{equation*}
Therefore, $\sin(e^{0.14})=0.912875$ which is correct upto the sixth place of decimal.
\section{conclusion}
Here we have defined a new type of calculus named G-calculus based on the idea of geometric differentiation defined by Grossman and Katz \cite{GrossmanKatz}. But Grossman and Katz took ordinary sum($+$) to produce increment to the independent variable $x$ such as $x_0, x_0+h, x_0+2h,...$ In that case some problem arise to discuss independently about the arithmetic system $(\oplus, \ominus, \odot, \oslash).$ That is why, idea of G-calculus comes our mind in which we took geometric sum($\oplus$) to produce increment to the independent variable $x$ such as $x_0, x_0\oplus h, x_0\oplus e^2\odot h,...$ Instead of mixing the ordinary arithmetic system($+,-,\times, \div $) and geometric arithmetic system $(\oplus, \ominus, \odot, \oslash)$, we are trying to formulate basic identities independently. As well as in \cite{KhirodBipan}, here, we are trying to bring up G-calculus to the attention of researchers in the different branches of analysis and its applications and advantages. We discussed in \cite{KhirodBipan} that the ordinary interpolation formulae are based upon the fundamental assumption that the data are expressible or can be expressed as a polynomial function with fair degree of accuracy. But geometric interpolation formulae have no such restriction. Because geometric interpolation formulae are based on geometric polynomials which are transcendental expressions in ordinary sense. So geometric interpolation formulae can be used to generate  transcendental functions, mainly to compute exponential and logarithmic functions.
\end{document}